\documentclass[a4paper,leqno]{article}

\usepackage{amsmath,amssymb,amscd,cite,authblk,multicol}

\allowdisplaybreaks

\begin{document}

\begin{center}

\textbf{\Large Projective (or spin) representations of finite groups. I}

\bigskip

{ Takeshi Hirai \footnote{Kyoto University, Yoshida-honmachi, Sakyo-ku, Kyoto 606-8501, JAPAN}~
\footnote{Hirai Mathematics Institute, 22-5 Nakazaichi-cho, Iwakura, Sakyo, Kyoto 606-0027, JAPAN}, Itsumi Mikami\footnote{Hirai Mathematics Institute}, Tatsuya Tsurii\footnote{Tokyo University of Information Sciences, 
Onaridai 4-1, Wakaba-ku, Chiba-shi,
Chiba, 265-8501, JAPAN},  and  Satoe Yamanaka\footnote{Department of Liberal Studies, 
National Institute of Technology, Nara College, 
22 Yata-cho, Yamatokoriyama, 
Nara 639-1080, JAPAN}}

\bigskip

\end{center}

\begin{abstract}
Schur multiplier $M(G)$ of a finite group $G$ has 
 been studied 
heavily.  To proceed further to the 
study of projective (or spin) representations of $G$ and their 
characters (called spin characters), 
it is necessary to construct explicitly a representation group $R(G)$ 
of $G$, a certain central extension of $G$ 
by $M(G)$, since projective representations of $G$ 
correspond bijectively to linear representations of $R(G)$. 
We propose here a practical method to construct $R(G)$ by 
 repetition of one-step efficient central extensions 
according to a certain choice of a series of 
elements of $M(G)$. This method 
 is also helpful for constructing linear 
representations of $R(G)$ and accordingly for 
calculating spin characters. 
Actually, we will apply this method to 
several examples of $G$ with prime number 3 in $M(G)$.

\bigskip

\medskip

\noindent
{\bf 2020 Mathematics Subject Classification:} Primary 20C25; Secondary 20F05, 20E99.

\noindent
{\bf Key Words:} efficient central extension of finite group, 
representation group, projective representation, 
spin representation, spin character.
\end{abstract}

\begin{multicols*}{2}

{\bf 1. Introduction. } 
Let $G$ be a finite group. A projective (or spin) representation $\pi$ 
of $G$ is a map $G\ni g \mapsto \pi(g) \in GL(n,{\boldsymbol C}),
\;n=\dim \pi,$ 
satisfying 
\begin{gather*}
\label{2023-11-04-1}
\pi(g)\pi(h)=r_{g,h}\pi(gh)\quad(g,h\in G,\;r_{g,h}
\in{\boldsymbol C}^\times),
\end{gather*}
where the function $r_{g,h}$ on $G\times G$ is called a factor set of 
$\pi$, 
and the set of equivalence classes of factor sets forms, as is known, 
 the cohomology $H^2(G,{\boldsymbol C}^\times)$, which is called 
{\it Schur multiplier} of $G$ and denoted by $M(G)$. 
Every projective representation $\pi$ is induced naturally from a linear 
representation of a finite central 
extension of $G$. According to Schur \cite{Sch1}, there exists a 
certain finite central extension of $G$ which has such a property 
for all the projective representations. 
The one with the minimum order is called {\it Representation group} 
of $G$ and 
 denoted by $R(G)$. $R(G)$ is a certain central extension of $G$ 
by $M(G)$ and there may exist a finite number 
of non-isomorphic $R(G)$'s, but we choose any one fixed. 

For an irreducible projective representation $\pi$ of $G$, 
there corresponds 
a one-dimensional 
character $\chi$ of $M(G)$ uniquely in such a way that, 
lifting up $\pi$ as a linear representation $\widetilde{\pi}$ of $R(G)$, 
we have for any $z\in M(G)$,  
$$
\widetilde{\pi}(z)=\chi(z)I\quad(z\in M(G)),
$$
where $I$ denotes the identity operator. 
Then $\chi$ is called the {\it spin type} of $\pi$, 
and all the irreducible 
spin representations are classified according to their spin types.   

{\bf 2. Preliminaries.}
Let $A$ be an abelian group and 
\begin{gather}
\label{2023-11-06-1}
1 \longrightarrow A \longrightarrow H
\longrightarrow G\longrightarrow 1\quad ({\rm exact})
\end{gather}
be a short exact sequence such that $A$ is contained in the 
center 
$Z(H)$ of $H$. Then we call $H$ a central extension of $G$ by 
$A$. It is called {\it efficient} if 
 $A$ is contained in 
the commutator subgroup $[H,H]$ of $H$, i.e., 
\begin{gather}
\label{2023-11-06-2}
A\subset Z(H)\cap [H,H].
\end{gather}

\noindent 
In the case where $A:=\langle a \rangle$ is a cyclic group  
 generated by $a$,  
we call $H$ a {\it one-step 
 central extension} of $G$, given by $a$.

A characterization of a representation group is given in 
Schur \cite[\S 5]{Sch1} as follows\footnote{
In Karpilovsky \cite[\S 2.1]{Kar}, \lq\lq\,a group 
$G^*$ is said to be {\it covering group} of $G$ 
if it satisfies (i) $A \subset Z(G^*)\cap [G^*,G^*]$\,; 
(ii) $A\cong M(G)$\,; (iii) $G\cong G^*/A$.\rq\rq \; 
This means that $G^*$ is exactly equal to 
{\it representation group} 
$R(G)$ in our terminology and notation. Here and 
also in our works 
\cite{THirai}, \cite{HiHo1}, \cite{HiHo2}, we use 
the term 
\,{\it covering group} of $G$\, simply for  
{\it central extension} (and in case 
appropriate, the term {\it universal covering group} for 
{\it representation group}, even though it is not necessary 
unique).
}
: 
{\it a central extension $H$ in\, {\rm (\ref{2023-11-06-1})} 
is isomorphic to a representation group 
if it is efficient and $A\cong M(G)$, that is, 
}
\begin{gather}
\label{2023-11-07-1}
A\cong M(G)\subset Z(H)\cap [H,H],\;\; H/M(G)\cong G.
\end{gather}

We note here the next fact about commutator $[x,y]:=xyx^{-1}y^{-1}$, which 
will have an important  meaning just below:   
\vskip.9ex

{\bf Lemma 2.1.}\;(Lemma 2.9.2 in \cite{Kar}) 
 {\it Suppose that in a group $H$, the commutator 
$z=[x,y]$ commutes with both $x$ and $y$, then for any $m,n\ge 1$,
} 
$$
[x^m,y^n]=z^{mn}.
$$

{\bf 3. A practical method of constructing representation group. } 
At present, a general method of constructing representation 
group of a finite group 
$G$ is not available, except for the cases of abelian groups and 
of a certain type of semidirect product groups. 
Thus we propose here a practical method 
of constructing a representation group of
 $G$. 

Now let $H$ in (\ref{2023-11-06-1}) be an efficient 
 central extension of $G$. 
\vskip.9ex

{\bf Definition 3.1. Efficient one-step 
central extension through a commuting pair.} 

Let $x,y\in H$ be a 
pair commuting with each other,  i.e., $[x,y]=1$, 
and with the same order $d>1$ (cf. Lemma 2.1). We call  
 a central extension $K$ of $H$\,:
$$
1\rightarrow Z \rightarrow K \stackrel{\alpha}{\rightarrow} H 
\rightarrow 1, \quad Z=\langle z\rangle,  
$$ 
as a {\it one-step efficient central extension} given by $z$, 
if the following two conditions hold:

 (A) The covering map $\alpha$ maps 
$\xi\stackrel{\alpha}{\to} x,\;\eta\stackrel{\alpha}{\to} y,$ 
and    
\;$[\xi,\eta]=z\stackrel{\alpha}{\to} [x,y]=1$.   
 The order of $z$ is $d$.  
\vskip.9ex

(B) \hspace*{9ex}
$Z \subset Z(K)\cap [K,K].$ 

\vskip.9ex

{\bf The method of construction proposed.} 
Starting from $G$ itself, repeating one-step efficient 
 central extensions through commuting pairs, well step by step, 
we may arrive by chance to the situation 
(\ref{2023-11-07-1}), giving a representation 
group $R(G)$.
 \vskip.9ex

{\bf Explanation of the method.}\;  
This method, maybe called as \lq\lq\,Efficient stairway up 
to {\it the Sky}\,\rq\rq, 
 where {\it the Sky}\, means the world of 
representation groups $R(G)$ (and also of spin representations 
and spin characters). 
This upward steps to {\it the Sky}\, 
is actually a backward process of 
descending downward steps on the same stairway from $R(G)$ to $G$. 

 Let's give a simplified explanation: take a decomposition of 
$M(G)$ as a product of cyclic groups: 
$$
M(G)=Z_1\times Z_2\times\cdots\times Z_r,\quad
 Z_j=\langle z_j\rangle,
$$
and take descending chain 
$$
R(G)\to H_1\to H_2 \to \;\cdots\;\to\, H_r\cong G,
$$
with 
$H_i=R(G)/Z^{(i)},\;Z^{(i)}:=(\prod_{j=1}^i Z_j),$ and natural 
homomorphism 
$$
\alpha_i: H_i \to H_{i+1},\:z_{i+1}
\stackrel{\alpha_i}{\longrightarrow} 1.
$$
We see from\; $M(G)\subset [R(G),R(G)]$\; that 
 $z_{i+1}\in [H_{i+1},H_{i+1}]$. So the central 
element $z_{i+1}$ is a 
product of several number of 
commutators $[x_k,y_k]$ with $x_k,y_k\in H_{i+1}$. 
Nevertheless we may expect a good luck that $z_{i+1}$ is 
expressed by a single one.

Among other things, we can expect that each upward 
step $(H_{i+1}, z_{i+1}) \uparrow H_i$ (i.e., commuting pair efficient 
central extension in {\it Efficient stairway to the Sky})
 helps to construct a 
complete set of spin irreducible 
representations of $H_i$ from that of $H_{i+1}$ by means of classical 
method 
in the case of semidirect product \cite{THirai}, and thus repeating 
step by step this construction,  finally to  
obtain 
such a complete set of $R(G)$. If succeeds, in this 
simplest example for explanation, we have 
 $M(H_i)=\prod_{k>i}Z_k$\,.

This method may be compatible 
with computer-aided calculations, because the finding of an 
efficient stairway 
to {\it the Sky}\, needs a lot of trial and error, 
success not necessarily guaranteed cf. \cite{Tsurii}.
\vskip.9ex

 {\bf Application of the method.} 
Construction of representation group and of a complete list of its  
irreducible linear 
representations has been rather 
well developed in the case where 
the prime number 2 appears in Schur multiplier 
$M(G)$, cf. e.g. \cite{HiHo1}, \cite{HiHo2}. 
So, we proceed to the cases where other primes 
appear in $M(G)$. We take several examples inspired by a 
list in the last part of \cite{Taha} and 
reproduced in \cite[Table 8.1]{Kar}.

{\bf 4. Case of a non-abelian group of order 18. } 
Let $G:=G_{20}$ be a non-abelian group 
of order 18  in \cite[Table 8.1]{Kar}.
 Its presentation is as follows:

{\bf Set of generators:} \; $\{a,\;b,\;c\}$;\quad 

{\bf Set of fundamental relations:} 
\begin{eqnarray}
\label{2023-11-09-1} 
\left\{ 
\begin{array}{ll}
a^2=b^2=c^2=1, 
\\
(abc)^2=1, \quad 
(ab)^3=(ac)^3=1.
\end{array}
\right.
\end{eqnarray}

Its Schur multiplier is given as $M(G_{20})={\mathbb Z}_3$.

{\bf 4.1. Trial to follow Schur's method.}\; 

First we tried to follow Schur's method for 
the symmetric groups ${\mathfrak S}_n$ in his 
epoch making work \cite{Sch3}. 

For an irreducible projective representation $\pi$ of 
$G=G_{20}$, put 
\;$
{\boldsymbol A}:=\pi(a),\;{\boldsymbol B}:=\pi(b),\;
{\boldsymbol C}:=\pi(c).
$
\\
The fundamental relations are represented,with constant 
 factors $\alpha, \beta,\cdots,\mu$, as  
\begin{gather*}
{\boldsymbol A}^2=\alpha I,\;{\boldsymbol B}^2=\beta I,\;
{\boldsymbol C}^2=\gamma I, \; 
({\boldsymbol A}{\boldsymbol B}{\boldsymbol C})^2=\delta I,
\\
({\boldsymbol A}{\boldsymbol B})^3=\lambda I,\;
({\boldsymbol A}{\boldsymbol C})^3=\mu I,
\end{gather*}
where $I$ denotes the identity operator. 
Then we try to simplify this set of relations and 
obtain  

\vskip.9ex
{\bf Theorem 4.1.}\; (i)\; {\it 
Through replacements of 
 ${\boldsymbol A},{\boldsymbol B}, {\boldsymbol C}$ by   
their constant multiples, the above set of relations among 
them turns to the following 
simplified form: 
\begin{gather*}
{\boldsymbol A}^2= {\boldsymbol B}^2= 
{\boldsymbol C}^2= I,\; 
({\boldsymbol A}{\boldsymbol B}{\boldsymbol C})^2=\tau I,
\\
({\boldsymbol A}{\boldsymbol B})^3=
({\boldsymbol A}{\boldsymbol C})^3= I,
\end{gather*}
where $\tau^3=1$, and therefore $M(G_{20})={\mathbb Z}_3$. 
}

(ii)\; {\it Let ${\mathfrak G}$ be an algebraic system with 
presentation given as 
\\
\quad{\rm \bf Set of generators:}\quad$\{\eta_1,\eta_2,\eta_3,\zeta\}${\rm ;} 
\\
\quad{\rm \bf Set of fundamental relations:} 
\begin{gather}
\label{}
\quad
\left\{
\begin{array}{ll}
\eta_1^{\,2}=\eta_2^{\,2}=\eta_3^{\,2}=1, 
\\
 (\eta_1\eta_2\eta_3)^2=\zeta,\quad 
\zeta\;\,{\rm central}, \;\zeta^3=1, \;\;
\\
(\eta_1\eta_2)^3 =(\eta_1\eta_3)^3=1.
\end{array}
\right.
\end{gather}
Then,  ${\mathfrak G}$ is a representation group of $G_{20}$, 
denoted by $R'(G_{20})$,  with 
the covering map 
$$
\Psi: \eta_1\to a,\;\eta_2\to b,\;\eta_3\to c,\;\zeta\to 1.
$$

}

{\it Proof.} (i) This is proved by calculations. 

(ii) Almost automatic from (i). 
\hfill
$\Box$
\vskip.9ex

Put $\omega=\exp(2\pi i/3)$, then 
$\tau=\omega^k\:(k=-1,0,1)$, and call $\tau$ the 
{\it spin type} of the spin representation $\pi$. 

Now we ask how to construct the triplet 
${\boldsymbol A},{\boldsymbol B}, {\boldsymbol C}$.
 The case of $\tau=1$ is the non-spin case, and 
we can apply for instance Mackey's induced 
 representation method.

{\bf 4.2. Case of symmetric groups. }   
In \cite[p.162]{Sch3}, Schur takes 
the presentation of $\mathfrak{S}_n$ as 
\\
\quad{\bf Set of generators:} 

$\{\mbox{\rm simple permutations} 
\;s_i=(i\;\,i{+}1); 1\le i<n\}$;
\\
\quad{\bf Set of 
fundamental relations:} 
\begin{gather*}
s_i^{\,2}=1\;(1\le i <n),\;(s_is_{i+1})^3=1\;(1\le i<n{-}1),
\\
 (s_is_k)^2=1\quad(|i-k|>1).
\end{gather*}

{\bf Theorem 4.2.} (i) (Cf.\;\cite[p.163]{Sch3}) 
{\it 
For an irreducible projective representation $\pi$, 
its normalised form  for 
the operators ${\boldsymbol C}_i=\pi(s_i)$ is 
\begin{gather*}
\label{}
{\boldsymbol C}_i^2=I\;(1{\le} i{<}n),\;
({\boldsymbol C}_i{\boldsymbol C}_{i+1})^3=I\;(1{\le}i{<}n{-}1),
\\
({\boldsymbol C}_i{\boldsymbol C}_k)^2=jI\;(|i-k|>1),
\end{gather*}
with $j^2=1$, and therefore $M(\mathfrak{S}_n)={\mathbb Z}_2\cong \{\pm 1\}$.
}

(ii) (Cf.\;\cite[p.164, (II')]{Sch3}) 
{\it  
Let ${\mathfrak T}'$ be an algebraic system with 
presentation given as follows: 
\\
\quad{\bf Set of generators:}\quad$\{\eta_i\;(1\le i<n),\; \zeta\}${\rm ;} 
\\
\quad{\bf Set of fundamental relations:} 
\begin{gather}
\nonumber
\left\{
\begin{array}{ll}
\eta_i^{\,2}=1\;(1{\le} i{<}n),\;(\eta_i\eta_{i+1})^3=1\;(1{\le} i{<}n{-}1),
\\
\zeta\;\;{\rm central}, \;\zeta^2=1, 
\\
(\eta_i\eta_k)^2 =\zeta\quad(|i-k| > 1).
\end{array}
\right.
\end{gather}
Then,  ${\mathfrak T}'$ is a representation group of ${\mathfrak S}_n$ 
 with the covering map 
$$
\Psi: \eta_i\to s_i\;(1\le i< n),\;\;\zeta\to 1.
$$
}

To construct irreducible 
projective representations $\pi$, 
or to get 
a set of operators  
${\boldsymbol C}_i\;(1\le i<n)$, 
Schur invented 
the very basic one,  
called {\it Hauptdarstellung zweiter Art}, with which 
every irreducible ones can be constructed by means of 
standard techniques. It works also for other groups $G$ 
whose Schur multiplier contains ${\mathbb Z}_2$ 
(Cf. \cite{HiHo1}, \cite{HiHo2}).

\vskip.9ex

Now return to our present case of $G_{20}$ or of 
the set of ${\boldsymbol A},{\boldsymbol B}, {\boldsymbol C}$. 
For spin case of $\tau=\omega^{\pm 1}\ne 1$, we don't have 
any 
such basic representation as {\it Hauptdarstellung} as above, 
and we understand that 
it is necessary to find new {\it Hauptdarstellung}\,  
for typical cases where Schur multiplier contains the cyclic group 
$C_3$ of order 3 (or $C_p$ of another prime order $p$).

{\bf 4.3. New way to construct $R(G_{20})$. }  
As is explained above, we are forced to invent a new method 
to construct representation group, using determined data 
of Schur multiplier. 

Now, for the  present use of a newly proposed 
 {\it one-step excellent central 
extension} of $G=G_{20}$, we look for a good pair 
$\{x_1,x_2\}$ in it commuting with each other. 
Take $x_1=ab$ and $x_2=ca$, then 
\begin{gather*}
x_1x_2=abca,\;x_2x_1=cb,\;{\rm and}\;(abc)^2=1,
\\ 
\;\therefore\; 
[x_1,x_2]=1.
\end{gather*}
So, we 
 introduce another presentation of $G$ as 

{\bf Set of generators:} \; $\{x_1=ab,\;x_2=ca,\;a\}$;

{\bf Set of fundamental relations:} 
\begin{eqnarray}
\label{2023-11-09-2}
\left\{ 
\begin{array}{ll}
a^2=1, \;x_1^{\;3}=x_2^{\;3}=1, 
\\
(ax_1)^2=(x_2a)^2=1,\;x_1x_2=x_2x_1. 
\end{array}
\right.
\end{eqnarray}

The subgroup $A=\langle a\rangle$ acts on $G$ as
 $\varphi(a)g:=aga^{-1}\;(g\in G)$.

\vskip.9ex
{\bf Lemma 4.3.} {\it 
The presentation in {\rm (\ref{2023-11-09-1})} 
of the group $G=G_{20}$ 
in {\rm \cite[p.278]{Kar}}   
is equivalent to that in {\rm (\ref{2023-11-09-2})} above. 
The second line in 
{\rm (\ref{2023-11-09-2})} can be replaced by 
$$
\varphi(a)(x_1)=x_1^{\,-1},\;\varphi(a)(x_2)=x_2^{\,-1},\;
[x_1,x_2]=1.
$$

Put $X_i:=\langle x_i\rangle\;\;(i=1,2)$, 
 and denote by $C_k$ the cyclic group of 
order $k$. Then  
the structure of the group $G_{20}$ is given as follows and so 
its  GAPIdentity (i.e., 
Identity numbers for finite groups given by GAP system) 
of $R(G_{20})$ is {\rm [18,4]} and
}
 \;$G_{20}\cong G_{18}^{\;\,4}$ \;(Cf. \cite{WP}):
\begin{gather}
\label{2023-11-16-1}
G_{20}=(X_1\times X_2)\rtimes_\varphi A \cong 
(C_3\times C_3)\rtimes_\varphi C_2.
\end{gather}

{\it Proof.}\; The assertions are proved through calculations by hands. 
\hfill
$\Box$
\vskip.9ex

With this commuting pair $x_1,x_2$, we try to 
construct an efficient central extension 
$H\stackrel{\delta}{\rightarrow} G$ 
coming from $[x_1,x_2]=1$. Let $\xi_1, \xi_2, z_{12}$ and $w$ 
be elements of $H$ covering as 
\begin{gather*}
\xi_i\stackrel{\delta}{\to} x_i\;(i=1,2),\;\;w\stackrel{\delta}{\to}a,
\quad{\rm and}
\\
[\xi_1,\xi_2]=z_{12}\,\stackrel{\delta}{\to}\, [x_1,x_2]=1.
\end{gather*}

The set of all elements of $H$ is 
\begin{gather*}
\label{2023-11-10-1}
h=
z_{12}^{\;\,\beta}\xi_1^{\,\gamma_1}\xi_2^{\,\gamma_2}w^\sigma
=:(\beta,\gamma_1,\gamma_2,\sigma),
\\
0\le \beta,\gamma_1,\gamma_2\le 2\;({\rm mod}\;3),\;
0\le \sigma\le 1\;({\rm mod}\;2). 
\end{gather*}  

The set of fundamental relations is 
\begin{gather}
\label{}
\left\{
\begin{array}{ll}
\xi_1^{\;3}=\xi_2^{\;3}=1,\;z_{12}^{\;\,3}=1,\;w^2=1,
\\
z_{12}\;\;{\rm central},\;\;\xi_2\xi_1=z_{12}^{\;-1}\xi_1\xi_2,
\\
\varphi(w)\xi_i=\xi_i^{\,-1}\;(i=1,2). 
\end{array}
\right.
\end{gather}

{\bf Lemma 4.4.} {\it 
The product rule in $H$ is given as  
\begin{gather}
\label{2023-07-18-1}
\nonumber
\big(z_{12}^{\;\beta} \xi_1^{\,\gamma_1}\xi_2^{\,\gamma_2}w^\sigma\big)
\big(z_{12}^{\;\beta'} \xi_1^{\,{\gamma'_1}}\xi_2^{\,{\gamma'_2}}
w^{\sigma'}\big)
=z_{12}^{\;\beta''} \xi_1^{\,{\gamma''_1}}\xi_2^{\,{\gamma''_2}}
w^{\sigma''},
\\[1ex]
\label{2023-07-18-2}
\left\{
\begin{array}{ll}
\beta''\equiv\beta+\beta'-\gamma_2 \cdot(-1)^\sigma
\gamma'_1\;&({\rm mod}\;3),
\\
\gamma''_1\equiv \gamma_1+(-1)^\sigma\gamma'_1&({\rm mod}\;3),
\\
\gamma''_2\equiv\gamma_2+(-1)^\sigma\gamma'_2 &({\rm mod}\;3),
\\
\sigma''\equiv \sigma+\sigma' &({\rm mod}\;2).
\end{array}
\right.
\end{gather}

Under this product rule, $H$ becomes actually a group. 
}

\vskip.9ex

{\it Proof.}\; Under the above product rule, we can prove  
by explicit calculations,  
existence of the identity element, and  
that of the inverse, and the 
associative law.  For instance, for associative law, we 
should 
prove the equality for $h_1,h_2,h_3\in H$,
$$
(h_1h_2)h_3=h_1(h_2h_3).
$$
So, calculate both sides under the rule (\ref{2023-07-18-2}) 
 by means of parameters, and then compare them.
\hfill $\Box$

\vskip.9ex
{\bf Theorem 4.5.} {\it The cyclic subgroup 
$Z_{12}:=\langle z_{12}\rangle$ is equal to the center of $H$, 
and isomorphic to 
$M(G),\;G=G_{20}\cong G_{18}^{\;\,4}$.  
Since 
$$
M(G)\cong Z(H)\subset Z(H)\cap [H,H],\;H/Z(H)\cong G,
$$ 
the group $H$ 
is a representation group $R(G)$. 

Structure of $R(G)$ is given 
with $\Xi_i:=\langle \xi_i\rangle\;\;(i=1,2)$ and $W=\langle w\rangle$ as 
follows and so its GAPIdentity is 
{\rm [54,5]} or {\rm [54,8]}\,: 
$$
R(G_{20})=\big[(\Xi_1\times Z_{12})\rtimes \Xi_2\big]\rtimes W 
\cong [(C_3\times C_3)\rtimes C_3]\rtimes C_2.
$$
}

{\bf Note 4.6.} The difference among two isomorphism classes  
$G_{54}^{\;\,5}$ and $G_{54}^{\;\,8}$ can be found in \cite{WP}
as follows: \lq\lq\,the former does not have 3-dimensional 
irreducible representations, whereas the latter does.\rq\rq 

Furthermore, in the succeeding paper, we will 
determine complete 
sets of irreducible representations of each of 
$G_{54}^{\;\,5}$ and $G_{54}^{\;\,8}$.

\vskip.9ex

{\bf 5. Case of a non-abelian group of order 27. } 
The group $G=G_{39}$, given in \cite[Table 8.1]{Kar}, 
is presented as follows:  

{\bf Set of generators:}\qquad $\{a,\; b\}$; 

{\bf Set of fundamental relations:} 
\begin{gather}
\label{2023-11-10-11}
a^3=b^3=1,\;
(ab)^3=(a^{-1}b)^3=1.
\end{gather}
Its Schur multiplier is given as 
$M(G_{39})={\mathbb Z}_3\times {\mathbb Z}_3$.
\vskip.9ex

{\bf Lemma 5.1.} 
{\it 
Put $c=bab^{-1}a^{-1}=[b,a]$, then $c^3=1$ and 
$Z(G)=\langle c\rangle$. Moreover, $bab^{-1}$ commutes with $a$, 
and $aba^{-1}$ commutes with $b$.
}
\vskip.9ex

For our conveniences, we introduce new presentation as follows: 

{\bf Set of generators:}
$$
\{x_1=b,\;x_2=c,\;x_3=a\}; 
$$

{\bf Set of fundamental relations:}
\begin{eqnarray}
\label{2023-11-14-31}
\left\{
\begin{array}{ll}
x_1^{\,3}=x_2^{\,3}=x_3^{\,3}=1,
\\
x_2\;\;{\rm central},\;\;\;
[x_1,x_3]=x_2.
\end{array}
\right.
\end{eqnarray}

{\bf Lemma 5.2.} {\it 
The group $G_{39}$ has 3 cyclic subgroups as 
$X_i=\langle x_i\rangle$ of 
type $C_3$, and its structure is given as 
$$
G_{39}=(X_1\times X_2)\rtimes X_3\cong (C_3\times C_3)
\rtimes C_3\,,
$$ 
and so in GAPIdentity
 \;$G_{39}\cong 
G_{27}^{\;\,3}$, and 
$$
X_1 \;\stackrel{\rm commute}\longleftrightarrow\;X_2\;
\stackrel{\rm commute}\longleftrightarrow\;X_3\,.
$$
}

Now take the commuting pair $[x_1,x_2]=1$ and construct 
one-step efficient central extension 
$\widetilde{G}\stackrel{\delta}{\rightarrow}G$ relative to
 it as follows: 

{\bf Set of generators:} \quad $\{\xi_1,\xi_2,\xi_3,z_{12}\}$;

{\bf Covering map:}\; 
$$
\xi_i\stackrel{\delta}{\rightarrow}x_i\;(1 \le i\le 3),\;
[\xi_1,\xi_2]=z_{12}\stackrel{\delta}{\rightarrow}
[x_1,x_2]=1, 
$$

{\bf Set of fundamental relations:} 
\begin{gather}
\label{2023-11-11-11}
\left\{
\begin{array}{ll} 
\xi_1^{\;3}=\xi_2^{\;3}=
\xi_3^{\;3}=1, \;\;z_{12}^{\;\;3}=1,
\\
z_{12}=[\xi_1,\xi_2]\;\;{\rm central},
\\
\xi_2=[\xi_1, \xi_3],\quad[\xi_2,\xi_3]=1.
\end{array}
\right.
\end{gather}

{\bf Theorem 5.3.} {\it 
The algebraic system $\widetilde{G}$, 
with the set of elements 
\begin{gather*}
\label{2023-11-11-12}
h=z_{12}^{\;\,\beta}\,\xi_1^{\,\gamma_1}
\xi_2^{\,\gamma_2}\xi_3^{\,\gamma_3}=:
(\beta,\gamma_1,\gamma_2,\gamma_3),
\\
0\le \beta,\,\gamma_1,\,\gamma_2,\,\gamma_3\le 2\;\;({\rm mod}\;3),
\end{gather*}
and with the rule of product {\rm (\ref{2023-11-11-11})},  
 becomes a group of order 81. 
Put 
$$
Z_{12}:=\langle z_{12}\rangle,\;\;\Xi_j:=\langle \xi_j\rangle\;
(1\le j\le 3), 
$$
then the structure of $\widetilde{G}$ is given as follows 
and so in GAPIdentity $\widetilde{G}\cong G_{81}^{\,10}$\,:
\begin{gather*}
\label{2023-10-18-2}
Z(\widetilde{G})=Z_{12},\quad
\big[\widetilde{G},\widetilde{G}\big]
=\Xi_2\times Z_{12}, 
\\
\widetilde{G}=
\big[(Z_{12}\times \Xi_1)\rtimes \Xi_2\big]\rtimes \Xi_3\;\;
\\
\;\;\cong\big[(C_3\times C_3)\rtimes C_3\big]\rtimes C_3.
\end{gather*}

}

As seen in \S 4, 
the key point of the proof of this theorem 
is to establish the associative law 
for the product in $\widetilde{G}$.  
\vskip.9ex

{\bf  Proof of associative law.}\; 
We prove the associative law:  
\begin{gather}
\label{2023-11-14-2}
\quad
(h_1h_2)h_3 = h_1(h_2h_3)\quad(h_1,h_2,h_3\in \widetilde{G}),
\end{gather}
by two lemmas. 
To begin with, we prepare some notations. 
Put $\varphi(h)g:=hgh^{-1}\;(h,g\in\widetilde{G})$ and 
rewrite the third line of fundamental relations in
 (\ref{2023-11-11-11}) as follows: 
\begin{gather}
\label{2023-11-14-1}
\left\{
\begin{array}{ll}
\varphi(\xi_3)\xi_1 &= z_{12}\xi_1\xi_2^{\;2},
\\
\varphi(\xi_3)^2\xi_1 &= z_{12}^{\;2}\,\xi_1\xi_2,
\\
\varphi(\xi_3)\xi_1^{\;2} &=\xi_1^{\;2}\xi_2,
\\
\varphi(\xi_3)^2\xi_1^{\;2} &=\xi_1^{\;2}\xi_2^{\;2},
\end{array}
\right.
\end{gather}
and note that every right hand side does not 
contain $\xi_3$, that is, they all 
belong to the subgroup 
$K:=\langle \xi_1,\xi_2\rangle$ generated by $\{\xi_1,\xi_2\}$. 
The set of all elements of $\widetilde{G}$ is 
\begin{gather*}
\label{2023-11-14-11}
h=z_{12}^{\;\;\beta}\,
 \xi_1^{\,\gamma_1}\xi_2^{\,\gamma_2}\xi_3^{\,\gamma_3}
=:(\beta, \gamma_1,\gamma_2,\gamma_3),
\\
0\le \beta, \gamma_1,\gamma_2,\gamma_3\le 2\quad({\rm mod}\;3).
\end{gather*}

Without loss of generality, to prove the associative law, 
we can restrict ourselves to the cases where 
\begin{gather*}
h_1=(0,0,b_2,b_3),\;h_2=(0, b'_1,b'_2,b'_3),
\\ 
h_3=(0,b''_1,b''_2,0).
\end{gather*}

{\bf Lemma 5.4.} 
\begin{gather*}
\label{}
J:=h_1h_2=\xi_2^{\;b_2}\xi_3^{\;b_3}
(\xi_1^{\;b'_1}\xi_2^{\;b'_2}\xi_3^{\;b'_3})\;\;\;\;
\\
=\xi_2^{\;b_2}\cdot
\varphi(\xi_3^{\;b_3})\xi_1^{\;b'_1}\cdot 
\xi_3^{\;b_3}\xi_2^{\;b'_2}\xi_3^{\;b'_3}
\\
=
\xi_2^{\;b_2}\cdot
\varphi(\xi_3^{\;b_3})\xi_1^{\;b'_1}\cdot 
\xi_2^{\;b'_2}\xi_3^{\;b_3+b'_3}.
\\
L:=(h_1h_2)h_3=Jh_3\hspace*{24.5ex}
\\
=\big[\xi_2^{\;b_2}\cdot
\varphi(\xi_3^{\;b_3})\xi_1^{\;b'_1}\cdot 
\xi_2^{\;b'_2}\xi_3^{\;b_3+b'_3}\big]\cdot 
\xi_1^{\;b''_1}\xi_2^{\;b''_2}
\\
=\big[\xi_2^{\;b_2}\cdot
\varphi(\xi_3^{\;b_3})\xi_1^{\;b'_1}\cdot 
\xi_2^{\;b'_2}\big]
\cdot\varphi(\xi_3^{\;b_3+b'_3})\xi_1^{\;b''_1}\cdot 
\xi_2^{\;b''_2}\cdot \xi_3^{\;b_3+b'_3}.
\end{gather*}

{\bf Lemma 5.5.} 
\begin{gather*}
\label{}
Q:=h_2h_3=
(\xi_1^{\;b'_1}\xi_2^{\;b'_2}\xi_3^{\;b'_3})
(\xi_1^{\;b''_1}\xi_2^{\;b''_2}) \quad
\\
=\xi_1^{\;b'_1}\xi_2^{\;b'_2}\cdot\varphi(\xi_3)^{b'_3}
\xi_1^{\;b''_1}\cdot\xi_3^{\;b'_3}\xi_2^{\;b''_2}
\\
\;=
\xi_1^{\;b'_1}\xi_2^{\;b'_2}\cdot\varphi(\xi_3)^{b'_3}
\xi_1^{\;b''_1}\cdot 
\xi_2^{\;b''_2}\xi_3^{\;b'_3}.
\\
R:=h_1(h_2h_3)=h_1Q\hspace*{24.9ex}
\\
\;=
\xi_2^{\;b_2}\xi_3^{\;b_3}\cdot
\big[
\xi_1^{\;b'_1}\xi_2^{\;b'_2}\cdot\varphi(\xi_3)^{b'_3}
\xi_1^{\;b''_1}\cdot 
\xi_2^{\;b''_2}\xi_3^{\;b'_3}\big]
\\
=
\xi_2^{\;b_2}
\varphi(\xi_3)^{b_3}\xi_1^{\;b'_1}\cdot 
\xi_3^{\;b_3}\xi_2^{\;b'_2}\cdot\varphi(\xi_3)^{b'_3}
\xi_1^{\;b''_1}\cdot 
\xi_2^{\;b''_2}\xi_3^{\;b'_3}
\\
=
\xi_2^{\;b_2}
\varphi(\xi_3)^{b_3}\xi_1^{\;b'_1}\cdot 
\xi_2^{\;b'_2}\cdot\varphi(\xi_3)^{b_3+b'_3}
\xi_1^{\;b''_1}\cdot 
\xi_2^{\;b''_2}\xi_3^{\;b_3+b'_3}.
\end{gather*}

Proof of the associative law is completed by comparing the 
last terms in 
Lemmas 5.4 and 5.5, not counting $\xi_3^{\;b_3+b'_3}$ at the end.
Thus we can say that this law for 
$\widetilde{G}$ comes from the same law for the subgroup $K$. 
\hfill
$\Box$
\vskip.9ex

{\bf Note 5.6.}\; Suggested by our computational results  
 in \cite{Tsurii} by using 
GAP system, we find other one-step efficient 
covering groups of $G=G_{39}$, 
 with slightly different {\bf Set of fundamental relations} 
comparing to (\ref{2023-11-11-11}) 
 of $\widetilde{G}$. The relations 
$$
\xi_1^{\;3}=1,\;\xi_3^{\;3}=1
$$ 
in (\ref{2023-11-11-11}) is replaced by 
$$
\xi_1^{\;3}=z_{12}^{\:\:a},\;\xi_3^{\;3}=z_{12}^{\;\;b}
$$ 
for any 
arbitrary choice of $a,b$ such that $0\le a, b\le 2$. 
Thus, we know that 
the choice of one-step efficient central extensions needs not 
be unique.

\vskip.9ex

{\bf 6. Case of a non-abelian group  $\widetilde{G}$ of order 81. } 
A one-step efficient central extension $H$ of $\widetilde{G}$ 
relative to the commuting pair 
$\xi_2$ and $\xi_3$, is 
of order 243, and supposed to be a 
representation group of $G_{39}$. 
$H$ is presented as follows: 
\vskip.9ex

{\bf Set of generators:}\qquad
$\{\eta_1,\eta_2,\eta_3,z_{12},z_{23}\}$; 

{\bf Covering map:} \;$\delta: H \to \widetilde{G}$, 
\begin{gather*}
\label{}
\eta_j\stackrel{\delta}{\rightarrow} \xi_j\;(1\le j\le 3),\;\;
z_{12}\stackrel{\delta}{\rightarrow}z_{12},
\\
[\eta_2,\eta_3]=z_{23}\stackrel{\delta}{\rightarrow}
[\xi_2,\xi_3]=1,
\end{gather*}

{\bf Set of fundamental relations:}\; 
\begin{gather*}
\label{2023-11-17-11}
\left\{
\begin{array}{ll}
\eta_j^{\;3}=1\;(1\le j\le 3),\;&
\eta_2\eta_1=z_{12}^{\;-1}\eta_1\eta_2,
\\[1.3ex]
\varphi(\eta_3)\eta_1=z_{12}\,\eta_1\eta_2^{\;2},\;
&
\varphi(\eta_3)\eta_1^{\;2}=
\eta_1^{\;2}\eta_2,
\\
\varphi(\eta_3^{\;2})\eta_1= z_{12}^{\;2}z_{23}\,\eta_1\eta_2,\;
&
\varphi(\eta_3^{\;2})\eta_1^{\;2}
=z_{23}^{\;\;2}\,\eta_1^{\;2}\eta_2^{\;2},
\\[1.3ex]
\varphi(\eta_3)\eta_2=z_{23}^{\;2}\eta_2,&
\varphi(\eta_3)\eta_2^{\;2}=z_{23}\,\eta_2^{\;2},
\\ 
\varphi(\eta_3^{\;2})\eta_2=z_{23}\,\eta_2,&
\varphi(\eta_3^{\;2})\eta_2^{\;2}=z_{23}^{\;2}\,\eta_2^{\;2}.
\end{array}
\right.
\end{gather*}

{\bf Important notice 6.1.}\; In the above fundamental relations, 
in the right hand side of each formula for 
$$
\varphi(\eta_3^{\;a})\eta_j^{\;b}\quad(1\le a,b\le 2,\,j=1,2),
$$ 
there does not appear 
the variable $\eta_3$, that is, they all 
belong to the subgroup $K$ generated by 
$\{z_{12},\,z_{23},\,\eta_1,\,\eta_2\}$.
\vskip.9ex

Let $H$ be an algebraic system generated by 
 $\{z_{12}, z_{23},\eta_1,\eta_2,\eta_3\}$\,: 
 
{\bf Set of elements of $H$\,:} 
\begin{gather*}
h
=z_{12}^{\;\beta_1}z_{23}^{\;\beta_2}\,\eta_1^{\,\gamma_1}
\eta_2^{\,\gamma_2}\eta_3^{\,\gamma_3}
=:(\beta_1,\beta_2,\gamma_1,\gamma_2,\gamma_3),
\\
0\le \beta_1,\beta_2,\gamma_1,\gamma_2,\gamma_3\le 2\;\;({\rm mod}\;3).
\end{gather*}

{\bf Theorem 6.2.} \;{\it 
Under the above product rule, the pair of\, {\rm \bf 
Set of elements} and\, {\rm \bf Set of fundamental relations}, the 
algebraic system $H$ 
gives a group of order 243. With 
$Z_{23}:=\langle z_{23}\rangle$ and  
$Y_j:=\langle\eta_j\rangle\;(1\le j\le 3)$, 
there holds 
\begin{gather*}
Z(H)=Z_{12}\times Z_{23},\quad
[H,H]
=Z_{12}\times Z_{23}\times Y_2,\;  
\end{gather*}
and the structure of the group is given as 
\begin{gather*}
H=
\big\{[(Z_{12}\times Y_1)\rtimes Y_2]\times Z_{23}\big\}
\rtimes Y_3\;\,
\\
\;\,\cong\big\{[(C_3\times C_3)\rtimes C_3]\times C_3\big\}\rtimes C_3.
\end{gather*}

The group $H$ is a representation group $R(G_{39})$ of 
$G_{39}\cong G_{27}^{\;\,3}$ 
of order 27, 
and is isomorphic to $G_{243}^{\;\;3}$ with 
GAPIdentity {\rm [243,3]}.
}
 
\vskip.9ex  

The key point of the proof of Theorem 6.2 is 
again to prove 
the associative law for the product.
\vskip.9ex

{\bf Proof of associative law.}\; 
To prove 
$$
(h_1h_2)h_3=h_1(h_2h_3)
$$ 
for $h_j\in H$, it is 
sufficient to restrict to the following case:  
\begin{gather*}
h_1= (0,0,0,b_2,b_3),\;h_2=(0,0,b'_1,b'_2,b'_3),
\\
h_3=(0,0,b''_1,b''_2,0).
\end{gather*}

The calculations of both sides are carried out 
as in Lemmas 5.4 and 5.5, 
almost a word-for-word repetition, except that here factors 
$\varphi(\eta_3^{\;a})\eta_2^{\;b}$ and $z_{23}^{\;\;c}$ appear. 

Anyhow, 
taking into account Important notice 6.1 above, 
 the associative law for $H$ 
is guaranteed by the same law for the subgroup $K$.

\vskip.9ex

{\bf Corollary 6.3.} {\it 
Let $\widetilde{G}\cong G_{81}^{\,10}$ be 
the group of order 81 in 
Theorem 5.3. Then 
$$
M(\widetilde{G})\cong Z_{23}\cong {\mathbb Z}_3\,. 
$$

}

\noindent
{\large \bf e-mail addresses:}

\medskip

\hspace{-2mm}Takeshi Hirai: hirai.takeshi.24e@st.kyoto-u.ac.jp

\hspace{-2mm}Itsumi Mikami: kojirou@kcn.ne.jp

\hspace{-2mm}Tatsuya Tsurii: t3tsuri23@rsch.tuis.ac.jp

\hspace{-2mm}Satoe Yamanaka: yamanaka@libe.nara-k.ac.jp

\end{multicols*}

\end{document}